\documentclass[12pt]{article}

\usepackage{amsfonts,graphicx}

\newtheorem{lemma}{Lemma}

\title{Detecting rigid convexity of\\
bivariate polynomials}

\author{Didier Henrion$^{1,2}$}

\date{Draft of \today}

\begin{document}

\maketitle

\footnotetext[1]{LAAS-CNRS, University of Toulouse, France}
\footnotetext[2]{Faculty of Electrical Engineering,
Czech Technical University in Prague, Czech Republic}

\begin{abstract}
Given a polynomial $x \in {\mathbb R}^n \mapsto p(x)$ in $n=2$
variables, a symbolic-numerical algorithm is first
described for detecting whether the connected component
of the plane sublevel set ${\mathcal P} = \{x : p(x) \geq 0\}$
containing the origin is rigidly convex, or equivalently,
whether it has a linear matrix inequality (LMI) representation,
or equivalently, if polynomial $p(x)$ is hyperbolic with
respect to the origin.
The problem boils down to checking
whether a univariate polynomial matrix is
positive semidefinite, an optimization problem that can be solved with eigenvalue
decomposition. When the variety ${\mathcal C} = \{x : p(x) = 0\}$
is an algebraic curve of genus zero, a second algorithm based
on B\'ezoutians is proposed to detect whether $\mathcal P$ has an LMI
representation and to build such a representation from a
rational parametrization of $\mathcal C$. Finally, some
extensions to positive
genus curves and to the case $n>2$ are mentioned.
\end{abstract}

\begin{center}
{\bf\small Keywords}\\[1em]
Polynomial, convexity, linear matrix inequality,
real algebraic geometry.\\
\end{center}

\section{Introduction}\label{intro}

Linear matrix inequalities (LMIs) are versatile modeling
objects in the context of convex programming, with
many engineering applications \cite{bn}.
An $n$-dimensional LMI set
is defined as
\begin{equation}\label{lmi}
{\mathcal F} = 
\{x \in {\mathbb R}^n \: :\: F(x) = F_0 + \sum_{i=1}^n x_i F_i
\succeq 0\}
\end{equation}
where the $F_i \in {\mathbb R}^{m\times m}$
are given symmetric matrices of size $m$
and $\succeq 0$ means positive semidefinite. From the
characteristic polynomial
\[
t \mapsto
\det(tI_m+F(x)) = p_0(x) + p_1(x)t + \cdots + p_{m-1}(x)t^m + t^m
\]
it follows from e.g. \cite[Theorem 20]{renegar} that
\begin{equation}\label{basic}
{\mathcal F} =
\{x \in {\mathbb R}^n \: :\: p_0(x) \geq 0, \ldots, p_{m-1}(x) \geq 0 \}.
\end{equation}
Hence the LMI set $\mathcal F$
is basic semialgebraic: it is the intersection
of polynomial sublevel sets. From
linearity of $F(x)$ and convexity of the cone of
positive semidefinite matrices, it also follows
that $\mathcal F$ is convex.
Hence LMI sets are {\em convex basic semialgebraic}.

\begin{figure}[h!]
\begin{center}
\includegraphics[width=15cm]{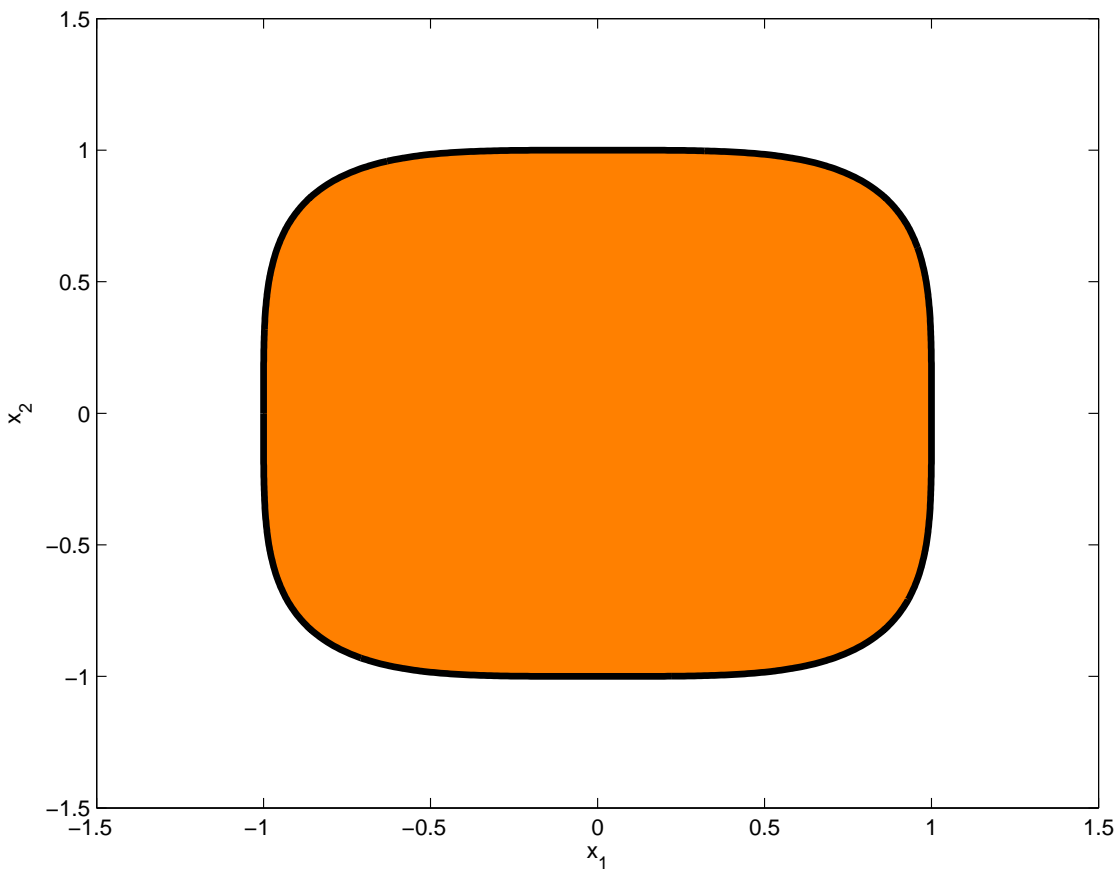}\\
\caption{The TV screen level set is not LMI.\label{tv_screen}}
\end{center}
\end{figure}

One may then wonder whether all convex basic semialgebraic
sets are LMI. In \cite{hv}, Helton and Vinnikov answer
by the negative, showing that in the plane ($n=2$)
some convex basic semialgebraic sets cannot be LMI.
An elementary example is the so-called TV screen set
defined by the Fermat quartic
\begin{equation}\label{tv}
\{x \in {\mathbb R}^2 \: :\: 1-x_1^4-x_2^4 \geq 0\}
\end{equation}
see Figure \ref{tv_screen}.

\subsection{Rigid convexity}

Assume that the set $\mathcal F$ defined in (\ref{lmi}) has a non-empty
interior, and choose a point $x_0$ in this interior, i.e.
\[
x_0 \in \mathrm{int}\:{\mathcal F} = \{x \: :\:
F(x) \succ 0\}
\]
where $\succ 0$ means positive definite.
A segment starting from $x_0$ attains the boundary of $\mathcal F$
when the determinant $p_0(x) = \mathrm{det}\:F(x)$
vanishes. The remaining polynomial inequalities
$p_i(x) \geq 0$, $i>0$ only isolate the convex connected component
containing $x_0$.
This motivated Helton and Vinnikov \cite{hv} to study
semialgebraic sets defined
by a single polynomial inequality
\begin{equation}\label{poly}
{\mathcal P} = \{x \in {\mathbb R}^n \: :\: p(x) \geq 0\}.
\end{equation}
The set $\{x \: :\: p(x) > 0\}$ is called an algebraic interior
with defining polynomial $p(x)$, and it is equal to $\mathrm{int}\:{\mathcal P}$
when $\mathcal P$ is convex. With these notations,
the question addressed in \cite{hv}
is as follows: what are the conditions satisfied
by a polynomial $p(x)$ so that $\mathcal P$
is an LMI set ?

For notational simplicity we will assume, without loss of generality,
that $x_0=0$, so that $\mathcal P$ contains the origin,
and hence we can normalize $p(x)$ so that $p(0)=1$.

If $p(x) = \mathrm{det}\:F(x)$ for
some matrix mapping $F(x)$ we say that $p(x)$ has
a determinantal representation. In particular,
the polynomial $p_0(x)$ in (\ref{basic}) has
a symmetric linear determinantal representation.

Consider an LMI set $\mathcal F$ as in (\ref{lmi}) and define
\[
p(x) = \det F(x)
\]
as the determinant of the symmetric pencil $F(x)$.
Note that $\deg p = m$, the dimension of $F(x)$.
Define the algebraic variety
\begin{equation}\label{var}
{\mathcal C} = \{x \in {\mathbb R}^n \: :\: p(x) = 0\}
\end{equation}
and notice that the boundary of $\mathcal F$ is
included in $\mathcal C$. Indeed, a point $x^*$
along the boundary of $\mathcal F$ is such that
the rank of $F(x^*)$ vanishes. Since the origin
belongs to $\mathcal F$ it holds $F_0 \succeq 0$.

Now consider a line passing through the origin,
parametrized as $x(t,z) = t z$ where $t \in \mathbb R$
is a parameter and $z \in {\mathbb R}^n$ is any
vector with unit norm. For all $z$, the symmetric matrix
$F(x(t,z)) = F_0 + t(z_1 F_1 + \cdots + z_n F_n)$
has only real eigenvalues as a pencil of $t$, and its
determinant $t \mapsto p(x(t,z)) = \det F(x(t))$ has only
real roots. Therefore, a given polynomial level set
$\mathcal P$ as in (\ref{poly}) is LMI only if the
polynomial $t \mapsto p(x(t,z))$ has only real roots for all $z$,
it must satisfy the so-called {\em real zero condition} \cite{hv}.
Geometrically it means that a generic line passing through
the origin must intersect the variety (\ref{var})
at $m = \deg p$ real points. The set $\mathcal P$ 
is then called {\em rigidly convex}, a geometric property
implying convexity.

A striking result of \cite{hv} is that rigid convexity
is also a sufficient condition for a polynomial level set
to be an LMI set in the plane, i.e. when $n=2$. For example,
it can be checked easily that the TV screen set
(\ref{tv}) is not rigidly convex since a generic
line cuts the quartic curve only twice.

In the litterature on partial differential equations,
polynomials satisfying real zero condition are also
called hyperbolic polynomials, and the corresponding
LMI set is called the hyperbolicity cone, see \cite{renegar}
for a survey, and \cite{lpr} for connections between
real zero and hyperbolic polynomials.

In passing, note the fundamental distinction between
an LMI set (as defined above) and a semidefinite representable
set, as defined in \cite{nn,bn}. A semidefinite representable
set is the projection of an LMI set:
\[
{\mathcal F} =
\{x \in {\mathbb R}^n \: :\: \exists u \in {\mathbb R}^{n_u}
\: :\: F(x,u) = F_0 + \sum_{i=1}^n x_i F_i + \sum_{j=1}^{n_u} u_j G_j
 \succeq 0\}
\]
where the variables $u_j$, sometimes called liftings, are
instrumental to the construction of the set through
an extended pencil $F(x,u)$. Such a set is called a
lifted LMI set. It is convex semialgebraic, but in general
it is not basic. However, it can be expressed as a union of basic
semialgebraic sets. In the case of the TV screen set
(\ref{tv}) a lifted LMI representation follows from
the extended pencil
\[
F(x,u) = \left[\begin{array}{cccccc}
1+u_1 & u_2 \\
u_2 & 1-u_1 \\
& & 1 & x_1 \\
& & x_1 & u_1 \\
& & & & 1 & x_2 \\
& & & & x_2 & u_2
\end{array}\right]
\]
obtained by introducing two liftings.
It seems that the problem of knowing which convex
semialgebraic sets are semidefinite representable
is still mostly open, see \cite{lasserre,hn}
for recent developments.

\subsection{Determinantal representation}

Once rigid convexity of a plane set, or equivalently the
real zero property of its defining polynomial, is established,
the next step is constructing an LMI representation.
Algebraically, given a real zero bivariate polynomial $p(x_1,x_2)$
of degree $m$,
the problem consists in finding symmetric matrices $F_0$, $F_1$
and $F_2$ of dimension $m$ such that
\[
p(x_1,x_2) = \det (F_0 + F_1 x_1 + F_2 x_2)
\]
and $F_0 \succeq 0$. If the $F_i$ are symmetric complex-valued
matrices, this is a well-studied problem
of algebraic geometry called {\em determinantal representation},
see \cite{room} for a classical reference
and \cite{beauville,piontkowski} for more recent surveys
and extensions to trivariate polynomials.

If one relaxes the dimension constraint (allowing the $F_i$
to have dimension larger than $m$) and the symmetry constraint
(allowing the $F_i$ to be non-symmetric), then results from
linear systems state-space realization theory (in particular 
linear fractional representations, LFRs) can be invoked to design
computer algorithms solving constructively
the determinantal representation problem. For example,
the LFR toolbox for Matlab \cite{lfr} is a user-friendly
package allowing to find non-symmetric determinantal
representations:
\begin{verbatim}
>> lfrs x1 x2
>> f=1/(1-x1^4-x2^4)
..
LFR-object with 1 output(s), 1 input(s) and 0 state(s).
Uncertainty blocks (globally (8 x 8)):
 Name  Dims  Type   Real/Cplx   Full/Scal   Bounds
 x1    4x4   LTI       r           s        [-1,1]
 x2    4x4   LTI       r           s        [-1,1]
\end{verbatim}
The software builds a state-space realization
of order 8 of the transfer function $f(x)=1/p(x)$.
This indicates that a non-symmetric real pencil $F(x)$
of dimension 8 could be found that satisfies
$\det F(x)=p(x)$, as evidenced by the following script
using the Symbolic Math Toolbox
(Matlab gateway to Maple):
\begin{verbatim}
>> syms x1 x2
>> D=diag([ones(1,4)*x1 ones(1,4)*x2]);
>> F=eye(8)-F.a*D
F =
[   1, -x1,   0,   0,   0,   0,   0,   0]
[   0,   1, -x1,   0,   0,   0,   0,   0]
[   0,   0,   1, -x1,   0,   0,   0,   0]
[ -x1,   0,   0,   1, -x2,   0,   0,   0]
[   0,   0,   0,   0,   1, -x2,   0,   0]
[   0,   0,   0,   0,   0,   1, -x2,   0]
[   0,   0,   0,   0,   0,   0,   1, -x2]
[ -x1,   0,   0,   0, -x2,   0,   0,   1]
>> det(F)
ans =
-x2^4+1-x1^4
\end{verbatim}
Note that LFR and state-space realization techniques
are not restricted to the bivariate case, but they
result in pencils of large dimension (typically much larger
than the degree of the polynomial), and there is apparently
no easy way to reduce the size of a pencil.

If one insists on having the $F_i$ symmetric, then
results from non-commutative state-space realizations
can be invoked to derive a determinantal representation,
at the price of relaxing the sign constraint on $F_0$.
An implementation is available in the {\tt NC}
Mathematica package \cite{hmv}. Here too, these techniques
may produce pencils of large dimension.

Now if one insists on having symmetric $F_i$ of minimal
dimension $m$, then two essentially equivalent constructive
procedures are known in the bivariate case to derive Hermitian
complex-valued $F_i$ from a defining polynomial $p(x_1,x_2)$
of degree $m$. Real symmetric solutions must then be extracted
from the set of complex Hermitian solutions.

The first one is based on the construction of a basis
for the Riemann-Roch space of complete linear systems 
of the algebraic plane curve $\mathcal C$ given
in (\ref{var}). The procedure is described in \cite{dixon}:
one needs to find a curve of degree $m-1$ touching $\mathcal C$
at each intersection point, i.e. the gradients must match.
The algorithm is illustrated in \cite{meyerbrandis}.
It is not clear however how to build a touching curve
ensuring $F_0 \succeq 0$. 

The second determinantal representation algorithm is
sketched in \cite{hv} and in much more detail in \cite{vinnikov93}.
It is based on complex Riemann surface theory \cite{griffiths,
fk}. Explicit expressions for the $F_i$ matrices are given via theta
functions. Numerically, the key ingredient is the computation
of the period matrix of the algebraic curve and the
corresponding Abel-Jacobi map. The period matrix of a curve
can be computed numerically with the {\tt algcurves} package of
Maple, see \cite{dvh} and the tutorial
\cite{deconinck} for recent developments, including new
algorithms for explicit computations of the Abel-Jacobi map.
A working computer implementation taking $p(x_1,x_2)$ as input
and producing the $F_i$ matrices as output is still missing
however.

\subsection{Contribution}

The focus of this paper is mostly on computational
methods and numerical algorithms. The contribution is twofold. 

First in Section \ref{rigid} we describe an algorithm for
detecting rigid convexity in the plane. Given a bivariate
polynomial $p(x_1,x_2)$, the algorithm uses a hybrid
symbolic-numerical method to detect whether the
connected component of the sublevel set (\ref{poly})
containing the origin is rigidly convex. 
The problem boils down to deciding whether a univariate
polynomial matrix is positive semidefinite.
This is a well-known problem in
linear systems theory, for which numerical linear
algebra algorithms are available (namely eigenvalue
decomposition), as well as a (more expensive but
more flexible) semidefinite programming formulation.

Then in Section \ref{rational} we describe an algorithm
for solving the determinantal
representation problem for algebraic plane curves of genus zero.
The algorithm is essentially symbolic, using B\'ezoutians,
but it assumes that a rational parametrization of the
curve is available. The idea behind the algorithm is
not new, and can be traced back to \cite{k},
as surveyed recently in \cite{kaplan}. An algorithm for
detecting rigid convexity of a connected component delimited by
such curves readily follows.

Extensions to positive genus algebraic plane curves and
higher dimensional sets are mentioned in Section \ref{extensions}.
In particular we survey the case of cubic plane curves and
cubic surfaces which are well understood. The case of quartic
(and higher degree) curves seems to be mostly open, and
computer implementations of determinantal representations
are still missing. Similarly, checking rigid convexity in
higher dimensions seems to be computational challenging
since it amounts to deciding whether a multivariate
polynomial matrix is positive semidefinite.

\section{Detecting rigid convexity in the plane}\label{rigid}

In this section we design an algorithm to assess whether the
connected component delimited by a bivariate polynomial around
the origin is rigidly convex. The idea is elementary and
consists in formulating algebraically the geometric condition
of rigid convexity of the set $\mathcal P$ defined in (\ref{poly}):
a line passing through the origin cuts
the algebraic curve $\mathcal C$ defined in (\ref{var}) a number of times which
is equal to the total degree $m$ of the defining bivariate polynomial
\[
x \in {\mathbb R}^2 \mapsto p(x) = \sum_{\alpha \in {\mathbb N}^2, |\alpha| \leq m}
p_{\alpha} x^{\alpha} = 
p_{00} + p_{10} x_1 + p_{01} x_2 +
p_{20} x_1^2 + p_{11} x_1 x_2 + \cdots
\]

A line passing through the origin can be
parametrized as:
\begin{equation}\label{trigo}
\begin{array}{rcl}
x_1 & = & r \cos \theta = t^{-1}(z^{-1}+z)\\
x_2 & = & r \sin \theta = it^{-1}(z^{-1}-z)
\end{array}
\end{equation}
where $z = e^{i\theta}$ and $t=2r^{-1}$.
Along this line, we define
\begin{equation}\label{pt}
t \in {\mathbb R} \mapsto q(t) = t^m p(x) = \sum_{k=0}^m q_k(z)t^k
\end{equation}
as a univariate polynomial of degree $m$
which vanishes on $\mathcal C$. Moreover $q(t)$ is monic
since $q_m(z) = p(0) = 1$. The remaining coefficients
are Laurent polynomials
\[
q_{\beta}(z) = \sum_{k=0}^m q_{\beta k}(z^k+z^{-k})
\]
with real coefficients, also called trigonometric
cosine polynomials.
Set $\mathcal P$ is rigidly convex if and only if this
polynomial has only real roots, i.e. if the number of
intersections of the line with the curve $\mathcal C$ is
maximal.

\subsection{Counting the real roots of a polynomial}

A well-known result of real algebraic geometry
\cite[Theorem 4.57]{bpr} states that a univariate polynomial
$q(t)$ of degree $m$ has only real roots if and only if
its Hermite matrix is positive semidefinite. The Hermite matrix
is the $m$-by-$m$ moment matrix of a discrete measure
supported with unit weights on the roots $x_1,\ldots,x_m$
of polynomial $q(t)$ (note that these roots are not necessarily
distinct). It is a symmetric Hankel matrix whose entries $(i,j)$
are Newton sums $N_{i+j} = \sum_{k=1}^m x_k^{i+j}$. The Newton
sums are elementary
symmetric functions of the roots that can be expressed
explicitly as polynomial functions of the coefficients of $q(t)$.
Recursive expressions are available to compute the $N_k$,
or equivalently, $N_k = \mathrm{trace}\:C^k$ where $C$
is a companion matrix of polynomial $q(t)$, i.e. a matrix
with eigenvalues $x_i$, see e.g. \cite[Proposition 4.54]{bpr}.
Recall that coefficients of the polynomial $q(t)$ given in (\ref{pt})
are Laurent polynomials. It follows that the Hermite matrix
of $q(t)$ is a symmetric trigonometric polynomial matrix of 
dimension $m$, that we denote by $H(z)$. We have proved
the following result.

\begin{lemma}\label{rh}
The bivariate polynomial $p(x)$ is rigidly convex
if and only if its Hermite matrix $H(z)$ is positive semidefinite
along the unit circle. Coefficients of $H(z)$
are explicit polynomial expressions of the coefficients of $p(x)$.
\end{lemma}

\subsection{Positive semidefiniteness of polynomial matrices}

The problem of checking positive semidefiniteness of a
polynomial matrix on the unit circle is generally referred
to as (discrete-time) spectral factorization. It is a well-known
problem of systems and circuit theory \cite{yakubovich,willems,glr}.
The positivity condition
can also be defined on the imaginary axis (continuous-time spectral
factorization) or the real axis. Various numerical methods are
available to solve this problem \cite{ks}. Several
algorithms are implemented in the Polynomial Toolbox for Matlab
\cite{polyx}. In increasing order of complexity, we can distinguish
between
\begin{itemize}
\item Newton-Raphson algorithms: the spectral factorization problem
is formulated as a quadratic polynomial matrix equation which is then solved
iteratively \cite{jk}. At each step, a linear polynomial matrix
equation must be solved \cite{hs}. Quadratic (resp. linear) convergence is
ensured locally if the polynomial matrix is positive definite (resp.
semidefinite);
\item polynomial operations: a sequence of elementary operations
is carried out in the ring of polynomials to reduce the
polynomial matrix to some canonical form, see \cite{callier}
and \cite{zh} for a recent survey. These algorithms are
cheap computationally but their numerical behavior (stability)
is unclear;
\item algebraic Riccati equation: using state-space realization,
the problem is formulated as a quadratic matrix equation,
which in turn can be solved via a matrix eigenvalue decomposition
with a particular structure \cite{willems,glr,tr};
\item semidefinite programming:
polynomial matrix positivity is formulated as a convex semidefinite
program, see \cite{tr} and the recent surveys \cite{g02,g03}.
The particular structure of this semidefinite program can be exploited
in interior-point schemes, in particular when forming
the gradient and Hessian. General purpose semidefinite
solvers can be used as well.
\end{itemize}

The semidefinite programming formulation of discrete-time
polynomial matrix factorization, a straightfoward transposition
of the continuous-time case studied in \cite{tr}, is as follows.
The symmetric trigonometric polynomial matrix
$H(z)=H_0+H_1(z+z^{-1})+\cdots+H_d(z^d+z^{-d})$
of size $m$
is positive semidefinite along the unit circle
if and only if there is a symmetric matrix $P$
of size $dm$ such that
\begin{equation}\label{spf}
\begin{array}{rcl}
L(P) & = &
\left[\begin{array}{c|ccc}
H_0 & H_1 & \cdots & H_d \\\hline
H_1 & 0 & & 0 \\
\vdots & & \ddots \\
H_d & 0 & \cdots & 0
\end{array}\right]
+
\left[\begin{array}{ccc}
I \\ \hline
& \ddots \\
& & I \\
0 & \cdots & 0
\end{array}\right]
P
\left[\begin{array}{c|ccc}
I & & & 0 \\
& \ddots & & \vdots \\
& & I & 0 
\end{array}\right] \\
&& -
\left[\begin{array}{ccc}
0 & \cdots & 0 \\ \hline
I \\
& \ddots \\
& & I \\
\end{array}\right]
P
\left[\begin{array}{c|ccc}
0 & I \\
\vdots & & \ddots\\
0 & & & I 
\end{array}\right] \\
& = & 
\left[\begin{array}{c|c}
H_0 & H_{01} \\ \hline
H_{01}^T & 0
\end{array}\right]
+
\left[\begin{array}{c}
B^T \\ A^T
\end{array}\right]
P
\left[\begin{array}{c|c}
B & A
\end{array}\right]
-
\left[\begin{array}{c}
D^T \\ C^T
\end{array}\right]
P
\left[\begin{array}{c|c}
D & C
\end{array}\right]
\succeq 0.
\end{array}
\end{equation}

Notice that the columns and rows of the above matrix
are indexed w.r.t. increasing powers of $z$
in such a way that
\[
B^T(z^{-1})L(P)B(z) = 
\left[\begin{array}{c}
I \\ z^{-1} \\ \cdots \\ z^{-d}
\end{array}\right]^T
L(P)
\left[\begin{array}{c}
I \\ z \\ \cdots \\ z^d
\end{array}\right] = H(z).
\]
Positive semidefiniteness of $L(P)$ then amounts
to the existence of a polynomial sum-of-squares
decomposition of $H(z)$. From the Schur decomposition
$L(P) = U^T U$ with $U = \left[\begin{array}{cccc}
U_0 & U_1 & \cdots & U_d \end{array}\right]$
it follows that
\begin{equation}\label{pspf}
H(z) = U(z^{-1})^T U(z).
\end{equation}
Polynomial matrix $U(z) = U_0 + U_1 z + \cdots + U_d z^d$
is called a spectral factor.

If the LMI problem (\ref{spf}) is feasible, then it
admits a whole family of solutions. Assuming that $H_0 \succ 0$,
maximizing the trace of $P$ subject to the
LMI constraint (\ref{spf}) yields a particular solution $P^*$
such that $\mathrm{rank}\:L(P^*) = m$.
It follows that the Schur complement of
\[
L(P) = 
\left[\begin{array}{c|c}
H_0+B^TPB-D^TPD & \star \\ \hline
H_{01}^T+A^TPB-C^TPD & A^TPA-C^TPC
\end{array}\right]
\]
w.r.t. $H_0$ vanishes, where symmetric entries are denoted by $\star$.
This means that $P^*$ satisfies the quadratic matrix equation
\[
A^TPA-C^TPC-(H_{01}^T+A^TPB-C^TPD)H_0^{-1}(H_{01}+B^TPA-D^TPC) = 0
\]
called the (discrete-time) algebraic Riccati equation.
In this case, the spectral factor $U(z)$ in (\ref{pspf})
is square non-singular.

\subsection{Example: cubic curve}

\begin{figure}[h!]
\begin{center}
\includegraphics[width=15cm]{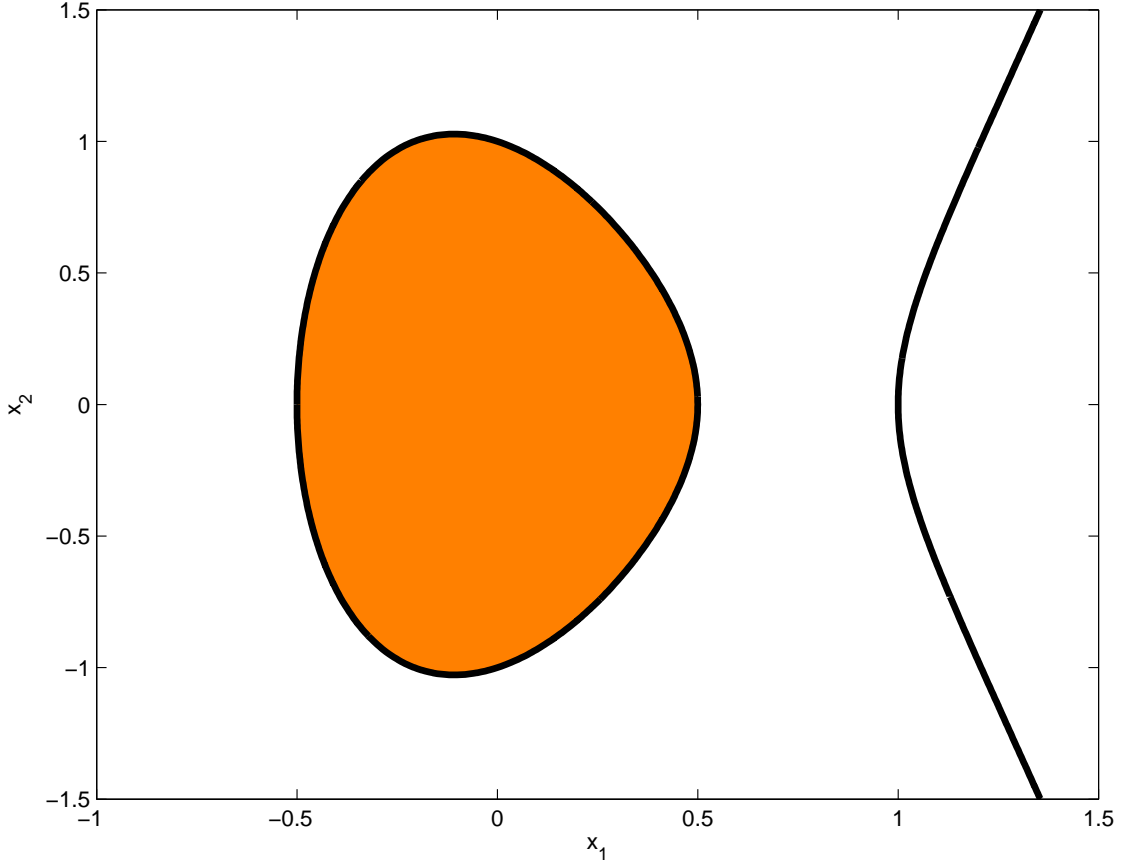}\\
\caption{Cubic curve and its component around the origin (shaded).\label{cubic_curve}}
\end{center}
\end{figure}

Consider the component of set (\ref{poly}) around the origin
delimited by the cubic polynomial $p(x)=1-x_1-4x_1^2-x_2^2+4x_1^3$,
see Figure \ref{cubic_curve}.

Using the substitution (\ref{trigo}) we obtain
\[
q(t) = 12(z+z^{-1})+4(z^3+z^{-3})-(10+3(z^2+z^{-2}))t-(z+z^{-1})t^2+t^3.
\]
From the companion matrix
\[
C = \left[\begin{array}{ccc}
z+z^{-1} & 10+3(z^2+z^{-2}) & -12(z+z^{-1})-4(z^3+z^{-3}) \\
1 & 0 & 0 \\
0 & 1 & 0
\end{array}\right]
\]
we build (symbolically) the Hermite matrix
\[
H(z) = \left[\begin{array}{ccc}
3 & \star & \star \\
z+z^{-1} & 22+7(z^2+z^{-2}) & \star \\
22+7(z^2+z^{-2}) & 6(z+z^{-1})-2(z^3+z^{-3}) &
250+124(z^2+z^{-2})+15(z^4+z^{-4})
\end{array}\right].
\]

Solving (numerically) the LMI (\ref{spf}) with SeDuMi interfaced with
YALMIP yields the spectral factorization
(\ref{pspf}) with factor (in Matlab notation)
\[
U(z) = \left[\begin{array}{ccc}
-0.9021-0.7094z^2 & -0.5284z+0.2027z^3 & -11.7639-9.6359z^2-1.5201z^4\\
0.1925z & 4.3449+1.6218z^2 & 0.7771z-0.5411z^3 \\
1.1578-0.5527z^2 & 0.3819z+0.1579z^3 & 2.4331-2.8689z^2-1.1844z^4
\end{array}\right]
\]
which certifies numerically that $p(x)$ is rigidly convex.

\subsection{Example: quartic curve}

Let us apply the algorithm to test rigid convexity
of the TV quartic level set (\ref{tv}) with
$p(x)=1-x_1^4-x_2^4$, see Figure \ref{tv_screen}.

We obtain $q(t) = -12-2(z^4+z^{-4})+t^4$ and
the Hermite matrix
\[
H(z) = \left[\begin{array}{cccc}
4 & \star & \star & \star \\
0 & 0 & \star & \star \\
0 & 0 & 48+8(z^4+z^{-4}) & \star \\
0 & 48+8(z^4+z^{-4}) & 0 & 0
\end{array}\right].
\]
From the zero diagonal entries and
the non-zero entries in the corresponding rows and columns
we conclude that $H(z)$ cannot be positive semidefinite
and hence that the TV quartic level set is not rigidly
convex.

\subsection{Numerical considerations}

Since the Hermite matrix $H(z)$ has a Hankel structure, and
positive definite symmetric Hankel matrices have a
conditioning (ratio of extreme eigenvalues) which can
be bounded below by an exponential function of the
matrix size \cite{becker,tyrty}, it may be appropriate
to apply a congruence transformation on matrix $H(z)$, also
called scaling.

For example, if $H(e^{i\theta_0})$ is positive definite
for some $\theta_0$ (say $\theta_0=0$, but other choices
are also possible), it admits
a Schur factorization $H(e^{i\theta_0})=V^T D V$ with $V$
orthogonal and $D$ diagonal non-singular. If $D$ is reasonably
well-conditioned, we can test positive semidefiniteness of
the modified trigonometric polynomial matrix
$H_0(z)=VD^{-1/2}H(z)D^{-1/2}V^T$ along the unit circle,
which is such that $H_0(e^{i\theta_0})$ is the identity matrix.
If $D$ is not well-conditioned, we can still use
$H_0(z) = VH(z)V^T$ which is such that $H_0(e^{i\theta_0})$ is
a diagonal matrix.

The impact of this data scaling
on the numerical behavior of the semidefinite programming or algebraic Riccati
equation solvers is however out of the scope of this paper.

\section{LMI sets and rational algebraic plane curves}\label{rational}

In the case that the algebraic curve $\mathcal C$ in (\ref{var})
has genus zero, i.e. the curve is rationally parametrizable,
an alternative algorithm can be devised to test rigid convexity
of a connected component delimited by $\mathcal C$. The algorithm
is based on elimination theory. It uses a particular symmetric
form of a resultant called the B\'ezoutian. As a by-product, the
algorithm also solves the determinantal representation problem
in this case. As surveyed recently in \cite{kaplan}, the key idea of
using B\'ezoutians in the context of determinantal
representations can be traced back to \cite{k}.

Starting from the implicit representation
\begin{equation}\label{implicit}
{\mathcal C} = \{x \in {\mathbb R}^2 \: :\: p(x) = 0\}
\end{equation}
of curve $\mathcal C$, with $p(x)$ a bivariate polynomial
of degree $m$, we apply a parametrization algorithm to obtain
an explicit representation
\begin{equation}\label{param}
{\mathcal C} = \{x \in {\mathbb R}^2 \: :\: x_1 = q_1(u)/q_0(u), \:
x_2 = q_2(u)/q_0(u), \: u \in {\mathbb R}\}
\end{equation}
with $q_i(u)$ univariate polynomials of degree $m$.
Algorithms for parametrizing an implicit algebraic curve are
described in \cite{abhyankar,sendra,vanhoeij}. An implementation
by Mark van Hoeij
is available in the {\tt algcurves} package of Maple. The
coefficients of $q_i(u)$ are generally found
in an algebraic extension of small degree over the field
of coefficients of $p(x)$.

With the help of resultants, we can eliminate the variable
$u$ in parametrization (\ref{param}) and recover an implicit
equation (\ref{implicit}), see \cite[Section 3.3]{cox}.
To address this implicitization problem,
we make use of a particular resultant, the B\'ezoutian, see 
\cite[Section 5.1.2]{mourrain}.
Given two univariate polynomials $g, h$ of the same degree $m$
(if the degree is not the same, the smallest degree polynomial is considered
as a degree $m$ polynomial with zero leading coefficients)
build the following bivariate polynomial
\[
\frac{g(u)h(v)-g(v)h(u)}{u-v} = \sum_{k=0}^{m-1} \sum_{l=0}^{m-1}
b_{kl} u^k v^l
\]
called the B\'ezoutian of $g$ and $h$,
and the corresponding symmetric matrix $B(g,h)$
of size $m\times m$ with entries $b_{kl}$ bilinear in coefficients
of $g$ and $h$. As shown e.g. in
\cite[Section 5.1.2]{mourrain},
the determinant of the B\'ezoutian matrix is the resultant,
so we can use it to derive the
implicit equation (\ref{implicit}) of a curve
from the explicit equations (\ref{param}).

\begin{lemma}\label{pencil}
Given polynomials $q_0, q_1, q_2$ in (\ref{param}), a polynomial
$p$ in (\ref{implicit}) is given by
$p(x) = \det F(x)$ where
\begin{equation}\label{px}
\begin{array}{rcl}
F(x) & = & B(q_1,q_2)+x_1B(q_2,q_0)+x_2B(q_1,q_0) \\
& = & F_0+x_1F_1+x_2F_2
\end{array}
\end{equation}
is a symmetric pencil of size $m$.
\end{lemma}

{\bf Proof:} Rewrite the system of equations (\ref{param}) as
\[
\begin{array}{rcccl}
g_1(u) & = & q_1(u) - x_1 q_0(u) & = & 0 \\
g_2(u) & = & q_2(u) - x_2 q_0(u) & = & 0 \\
\end{array}
\]
and use the B\'ezoutian resultant
to eliminate indeterminate $u$ and obtain conditions
for a point $(x_1,x_2)$ to belong to the curve. The B\'ezoutian
matrix is $B(g_1,g_2) = B(q_1-x_1q_0,q_2-x_2q_0) =
B(q_1,q_2)+x_1B(q_2,q_0)+x_2B(q_1,q_0)$. Linearity in $x$
follows from bilinearity of the B\'ezoutian
and the common factor $q_0(u)$.$\Box$

Lemma \ref{pencil} provides an implicit equation of
curve (\ref{implicit}) in symmetric linear determinantal form.

\subsection{Detecting rigid convexity}

Once polynomial $p(x)$ is in symmetric linear determinantal
form as in Lemma \ref{pencil}, checking rigid convexity
of the connected component containing the origin $x=0$
amounts to testing positive definiteness of $F(0) = F_0 = B(q_1,q_2)$.

\begin{lemma}\label{definite}
The B\'ezoutian matrix $B(q_1,q_2)$ is positive semidefinite
if and only if polynomial $q_1(u)$ and $q_2(u)$
have only real roots that interlace.
\end{lemma}

{\bf Proof:} The signature (number of
positive eigenvalues minus number of negative eigenvalues)
of the B\'ezoutian of $q_1(u)$ and $q_2(u)$
is the Cauchy index of the rational function $q_1(u)/q_2(u)$,
the number of jumps of the function from $-\infty$ to $+\infty$
minus the number of jumps from $+\infty$ to $-\infty$, see
\cite[Definition 2.53]{bpr} or \cite[Theorem 9.4]{bpr}.
It is maximum when $B(q_1,q_2)$ is positive definite.
This occurs if and only if the roots
of $q_1(u)$ and $q_2(u)$ are all real and interlace.$\Box$

\begin{lemma}\label{rigidbez}
The connected component around the origin delimited
by curve (\ref{implicit}) is rigidly convex if and only if
$B(q_1,q_2) \succeq 0$.
\end{lemma}

{\bf Proof:} Since $F_0 \succeq 0$, the set admits the LMI
representation $\{x \in {\mathbb R}^2 \: :\: F(x) \succeq 0\}$,
which is equivalent to being rigidly convex. $\Box$

\subsection{Finding a rigidly convex component}

If the connected component around the origin is not
ridigly convex, it may happen that there is another
rigidly convex connected component elsewhere.
To find it, it suffices to determine a point $x^*$
such that $F(x^*) \succeq 0$.
This is equivalent to solving
a bivariate LMI problem.

We can apply primal-dual interior-point methods \cite{nn}
to solve this semidefinite programming problem, since
the function $f(x) = - \log p(x) = \log \det F(x)^{-1}$ is a strictly
convex self-concordant barrier for the interior of the LMI set. If the
LMI set is bounded, minimizing $f(x)$ yields the analytic
center of the set. If the LMI set
is empty, the dual semidefinite problem yields a
Farkas certificate of infeasibility. However, in the
bivariate case a point $x^*$ satisfying
$F(x^*) \succeq 0$ can be found more easily
with real algebraic geometry and univariate
polynomial root extraction.

A first approach consists in identifying the local
minimizers of function $f(x)$. They are
such that the gradient $g(x)$ of $p(x)$ vanishes, i.e.
they are such that 
\begin{equation}\label{grad}
g_i(x) = \frac{\partial p(x)}{\partial x_i} =
\mathrm{trace}(p(x)F^{-1}(x)F_i) = 0.
\end{equation}
We can characterize these minimizers by eliminating
one variable, say $x_1$, from the system $g_1(x)=g_2(x)=0$,
and solving for the other variable $x_2$ via
polynomial root extraction. Resultants
can be used for that purpose.

A second approach consists in finding points on
the boundary of the LMI set, which are such that
$p(x)=0$ and either $g_1(x)=0$ or $g_2(x)=0$.
Here too, resultants can be applied to end up
with a polynomial root extraction problem.

From the points generated by these two procedures,
we keep only those satisfying $F(x) \succeq 0$,
an inequality that can be certified by testing
the signs of the coefficients of the characteristic
polynomial of $F(x)$, as explained in the
introduction.

\subsection{Example: capricorn curve}

Let $p(x) = x_1^2(x_1^2+x_2^2)-2(x_1^2+x_2^2-x_2)^2$.
With the {\tt parametrization} function of the {\tt algcurves}
package of Maple, we obtain a rational parametrization
\[
\begin{array}{rcl}
q_0(t) & = & 45-8t+10t^2+t^4 \\
q_1(t) & = & -7+44t-18t^2-4t^3+t^4\\
q_2(t) & = & 49-28t-10t^2+4t^3+t^4.
\end{array}
\]
With the {\tt BezoutMatrix} function of the {\tt LinearAlgebra}
package, we build the corresponding symmetric pencil
\[\small
F(x) = \left[\begin{array}{cccc}
8-4x_1-4x_2 & \star & \star & \star \\
8+20x_1-28x_2 & 40+60x_1+92x_2 & \star & \star \\
-72+20x_1+52x_2 & -8-36x_1-84x_2 & 776+540x_1+476x_2 & \star \\
56-4x_1-52x_2 & -168+180x_1+180x_2 & -952-940x_1+740x_2 & 1960-868x_1-1924x_2
\end{array}\right]
\]
whose determinant (up to a constant factor) is equal to $p(x)$.
The eigenvalues of $F(0)$ are equal to $0$ (double)
and $1392\pm48\sqrt{533}$. They are all non-negative
which indicates that the origin lies on the boundary of
an LMI region defined by $F(x) \succeq 0$.

Values of $x_2$ at local optima satisfying the system
of cubic equations (\ref{grad})
can be found with Maple as follows:
\begin{verbatim}
> p:=x1^2*(x1^2+x2^2)-2*(x1^2+x2^2-x2)^2:
> solve(resultant(diff(p,x1),diff(p,x2),x1));
0, 0, 0, 1, 1/2, 3+sqrt(5), 3-sqrt(5), 3+sqrt(5), 3-sqrt(5)
\end{verbatim}
from which it follows that, say, the point $x_1=0$, $x_2=1/2$
is such that $F(x) \succ 0$.

\begin{figure}[h!]
\begin{center}
\includegraphics{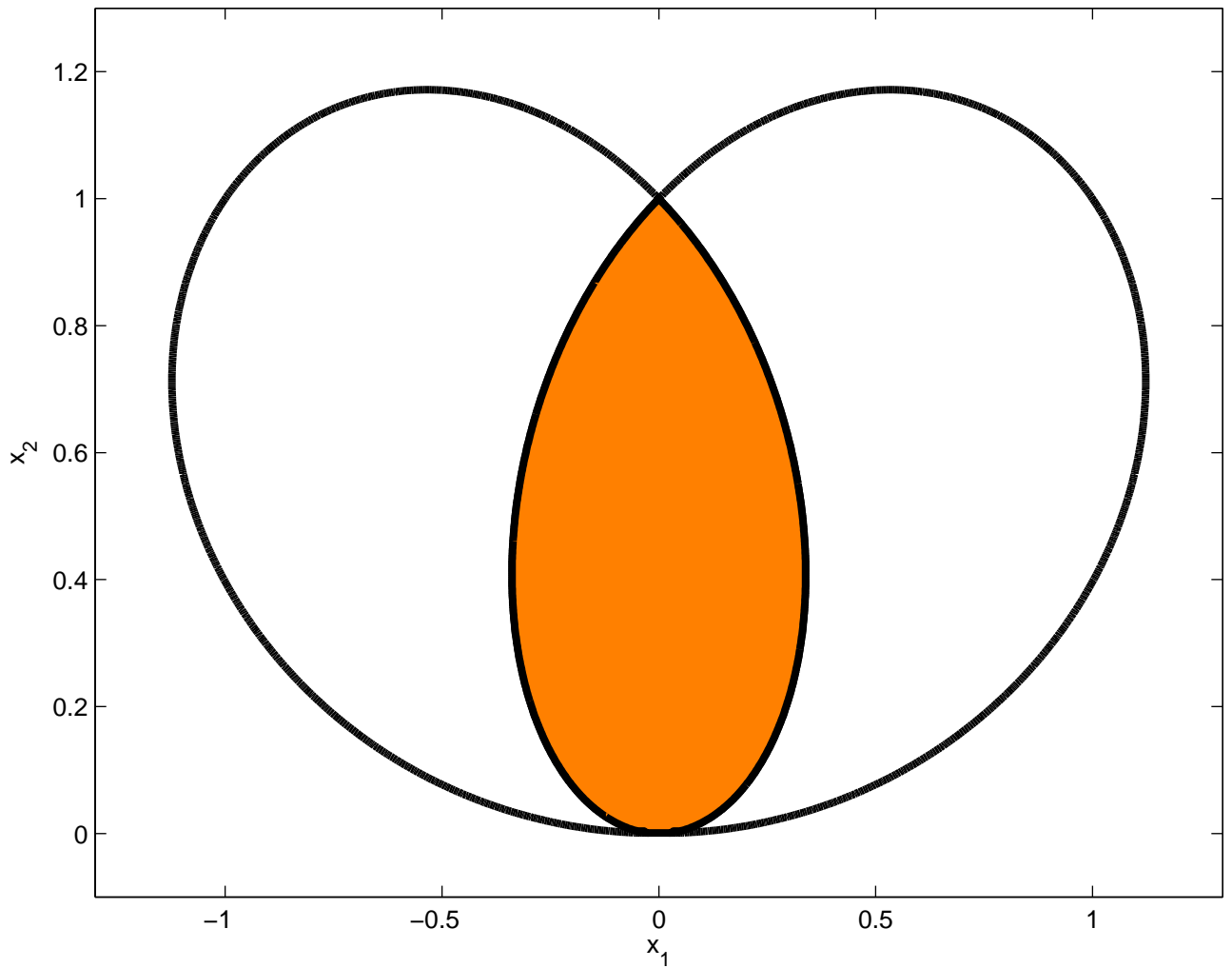}\\
\caption{Capricorn curve defining an LMI region (shaded).\label{capricorn_curve}}
\end{center}
\end{figure}

The corresponding LMI region together with the quartic
capricorn curve are represented on Figure
\ref{capricorn_curve}.

\subsection{Example: bean curve}

\begin{figure}[h!]
\begin{center}
\includegraphics[width=15cm]{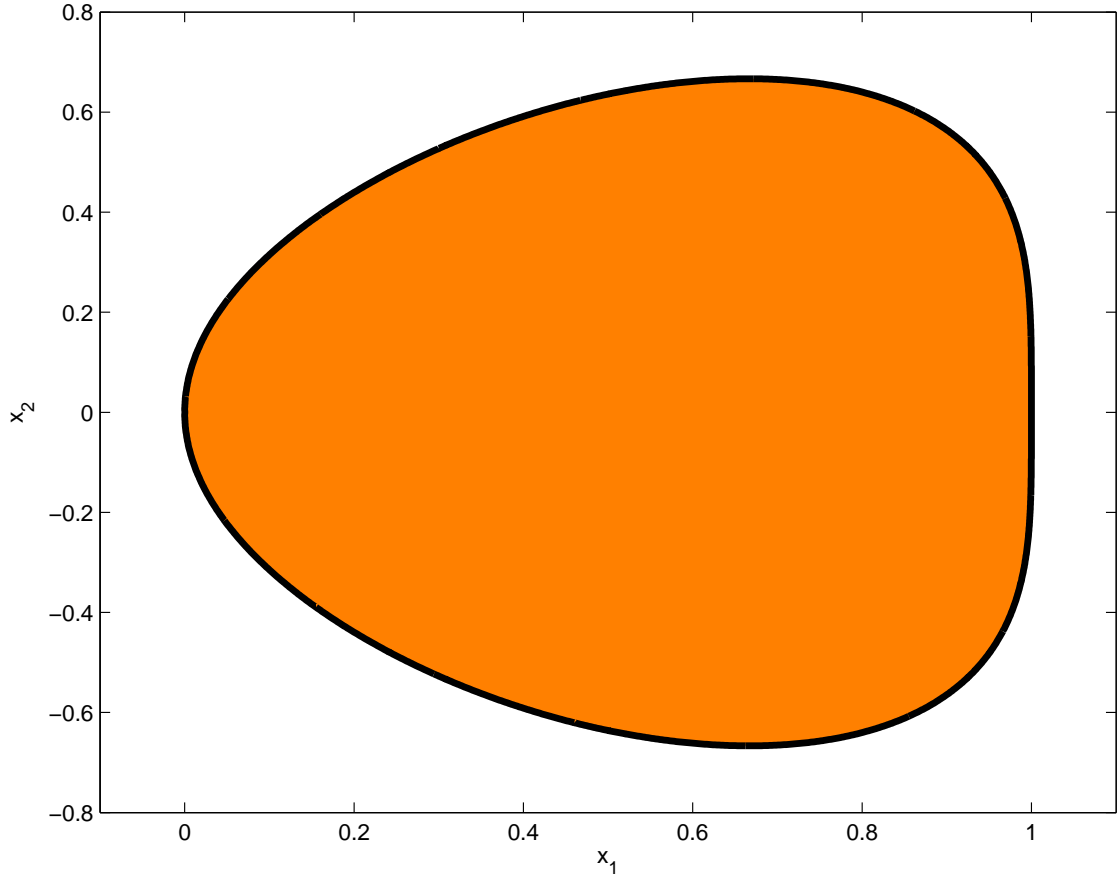}\\
\caption{Bean curve defining a region which is not LMI (shaded).
\label{bean_curve}}
\end{center}
\end{figure}

Let $p(x)=x_1^4+x_1^2x_2^2+x_2^4-x_1^3+x_1x_2^2$.
With the B\'ezoutians we obtain the following pencil
\[
F(x) = \left[\begin{array}{cccc}
x_1& \star & \star & \star\\
x_2&1& \star & \star \\
x_1&x_2&0& \star\\
x_2&1-x_1&0&1-x_1
\end {array}
\right].
\]
We can check that $F(0)$ has eigenvalues $2$ and $0$ (triple)
and this is the only point for which $F(x) \succeq 0$.
It follows that the convex set delimited by the curve
$p(x)=0$ is not LMI, see Figure \ref{bean_curve}.

\section{Extensions}\label{extensions}

In this paragraph we outline some potential extensions
of the results to algebraic plane curves of positive
genus and varieties of higher dimensions.

\subsection{Cubic plane curves}

The case of cubic plane algebraic curves is well understood,
see e.g. \cite{vinnikov86} or \cite{taussky}. Singular cubics
(genus zero) can be handled
via B\'ezoutians as in Section \ref{rational}.
Smooth cubics (genus one), also called elliptic curves,
can be handled via their Hessians.

Let $p(x_1,x_2)$ be a cubic polynomial that
we homogeneize to $p(x_0,x_1,x_2) = x_0^3p(x_1/x_0,x_2/x_0)$.
Define its 3-by-3 symmetric
Hessian matrix $H(p(x))$ with entries
\[
H_{ij} = \frac{\partial^2 p(x)}{\partial x_i \partial x_j}
\]
and the corresponding Hessian $h(x)=\det H(p(x))$.
The elliptic curve $p(x)=0$ has 9 inflection points,
or flexes, satisfying $p(x)=h(x)=0$, and 3 of them are real.
Since $p(x)$ and $h(x)$ share the same flexes and the Hessian
matrix yields a symmetric linear determinantal representation for $h(x)$,
we can use homotopy to find a determinantal representation
for $p(x)$.

For real $t$ define the parametrized Hessian
$g(x,t)=\det H(h(x)+tp(x))$ and find $t^*$ satisfying
$g(x^*,t^*)=p(x^*)$ at a real flex $x^*$ by solving
a cubic equation. As a result, we obtain three
distinct symmetric pencils not equivalent by
congruence transformation. One of the may be definite
hence LMI.

For example, let $p(x)=x_1^3-x_2^2-x_1$. Build the
Hessian $h(x)=\det H(p(x)) = 8(x_0^3+3x_0x_1^2-3x_1x_2^2)$
and the parametrized Hessian $g(x,t)=\det H(h(x)+tp(x))=
24t^3x_0x_1^2-576t^2x^2_0x_1+\cdots+110592x_1^3$.
Polynomial $g(x,t)$ matches $g(x)$ at flex $x^*_0=0$ for
$t^* \in \{0, 24, -24\}$
yielding the following three representations
\[
\begin{array}{rcl}
F^1(x) & = & \left[\begin{array}{ccc}
1 & \star & \star \\ -x_2 & -x_1 & \star \\ x_1 & 0 & 1
\end{array}\right] \\ \\
F^2(x) & = & 4^{-\frac{1}{3}}\left[\begin{array}{ccc}
1+3x_1 & \star & \star \\
-x_2 & -1-x_1 & \star\\
-1+x_1 & -x_2 & 1-x_1
\end{array}\right] \\ \\
F^3(x) & = & 4^{-\frac{1}{3}}\left[\begin{array}{ccc}
1-3x_1 & \star & \star\\
-x_2 & 1-x_1 & \star \\
1+x_1 & x_2 & 1+x_1
\end{array}\right]
\end{array}
\]
such that $\det F^i(x)=p(x)$ for all $i=1,2,3$.
Only the first one generates an LMI set $F^1(x)\succeq 0$.

\subsection{Positive genus plane curves}

The case of algebraic plane curves
of positive genus and degree equal to four (quartic)
or higher is mostly open. Whereas rigid convexity of higher
degree polynomials can be checked with the proposed approach,
there is no known implementation of an algorithm that
produces symmetric linear determinantal (and hence LMI)
representations in this case. For quartics, contact
curves can be recovered from bitangents. In \cite{edge}
complex symmetric linear determinantal representations
of the quartic $1+x^4_1+x^4_2$ could be derived from the
equations of the bitangents found previously
by Cayley for this particular curve.

B\'ezoutians can be generalized to the multivariate
case, as surveyed in \cite{mourrain}. In Lemma \ref{pencil}
we derived a symmetric linear determinantal representation
by eliminating the variable $u$ in the system of
equations
\[
\begin{array}{rclcl}
g_1(u) & = & q_1(u) - x_1 q_0(u) & = & 0 \\
g_2(u) & = & q_2(u) - x_2 q_0(u) & = & 0 \\
\end{array}
\]
corresponding to a rational parametrization
$x_1(u) = q_1(u)/q_0(u)$, $x_2(u) = q_2(u)/q_0(u)$
of the curve $p(x_1,x_2)=0$. In the positive genus case, such
a rational parametrization is not available, but
we can still define a system of equations
\[
\begin{array}{rclcl}
g_1(u_1,u_2) & = & x_1-u_1 & = & 0 \\
g_2(u_1,u_2) & = & x_2-u_2 & = & 0 \\
g_3(u_1,u_2) & = & p(u_1,u_2) & = & 0 
\end{array}
\]
describing the curve $p(x_1,x_2)=0$ after eliminating
variables $u_1$ and $u_2$. Define the discrete differentials
\[
\partial_1 g(u,v) = \frac{g(u_1,u_2)-g(v_1,v_2)}{u_1-v_1}, \quad
\partial_2 g(u,v) = \frac{g(v_1,u_2)-g(v_1,v_2)}{u_2-v_2}
\]
and the quadratic form
\[
\det \left[\begin{array}{ccc}
g_1 & \partial_1 g_1 & \partial_2 g_1 \\
g_2 & \partial_1 g_2 & \partial_2 g_2 \\
g_3 & \partial_1 g_3 & \partial_2 g_3
\end{array}\right] =
\det \left[\begin{array}{ccc}
x_1-u_1 & -1 & 0 \\
x_2-u_2 & 0 & -1 \\
p(u_1,u_2) & \partial_1 p(u,v) & \partial_2 p(u,v)
\end{array}\right] =
\sum_{\alpha,\beta} f_{\alpha,\beta} u^{\alpha}v^{\beta}
\]
using bi-indices $\alpha$ and $\beta$. Then the matrix
$F(x)$ of the quadratic form is a symmetric pencil
satisfying $\det F(x) = p(x)q(x)$ where $q(x)$
is an extraneous factor. We hope that $q(x)$ does
not depend on $x$, even though this cannot be
guaranteed in general. For example, in the case of the
Fermat curve $p(x) = 1-x^4_1-x^4_2$ whose genus is three,
using the {\tt multires} package for Maple \cite{multires},
we could obtain
\[
F(x) = \left[\begin{array}{ccccccc}
-1 & 0 & 0 & 0 & 0 & x_1 & x_2 \\
0 & 0 & 0 & x_1 & 0 & -1 & 0 \\
0 & 0 & 0 & 0 & x_2 & 0 & -1 \\
0 & x_1 & 0 & -1 & 0 & 0 & 0 \\
0 & 0 & x_2 & 0 & -1 & 0 & 0 \\
x_1 & -1 & 0 & 0 & 0 & 0 & 0 \\
x_2 & 0 & -1 & 0 & 0 & 0 & 0
\end{array}\right]
\]
which is such that $\det F(x) = -p(x)$, i.e. $q(x)=-1$.

\subsection{Surfaces and hypersurfaces}

The case $n=m=3$, i.e. cubic surfaces, is well
understood, see \cite{kosir} for a full constructive
development. All the self-adjoint linear determinantal
representations can be obtained from the tritangent
planes. The number of non-equivalent
representations depends on the number
and class of real lines among the 27 complex lines
of the surface. See \cite{szilagyi} for
a nice survey on cubic surfaces.

\begin{figure}[h!]
\begin{center}
\includegraphics[width=15cm]{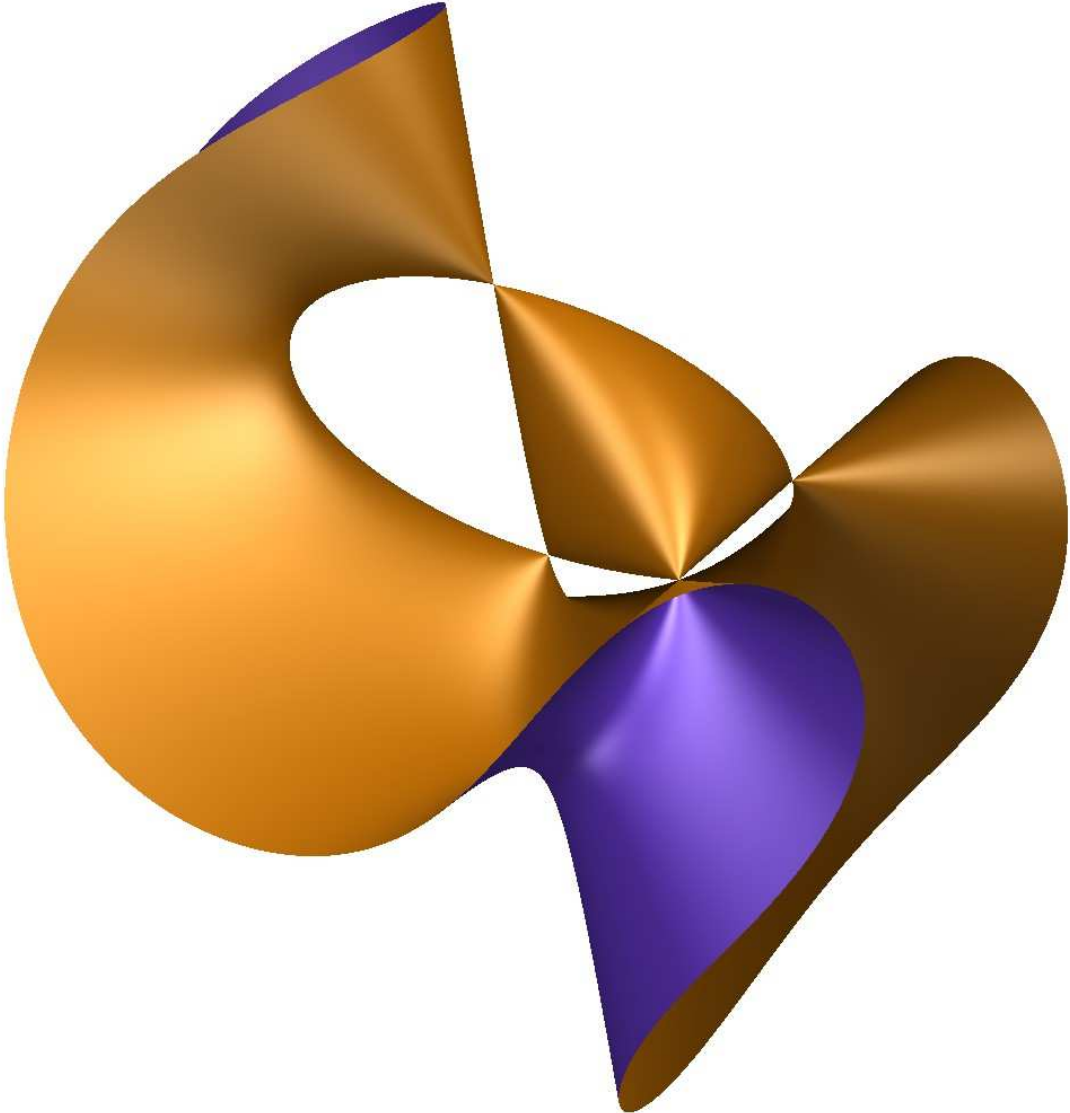}\\
\caption{The Cayley cubic surface with its convex connected
component.\label{cayley_cubic}}
\end{center}
\end{figure}

A well-known example is the Cayley cubic
\[
\frac{1}{u_0}+\frac{1}{u_1}+\frac{1}{u_2}+\frac{1}{u_3} = 0
\]
whose algebraic equation is
\[
u_0 u_1 u_2 + u_0 u_1 u_3 + u_0 u_2 u_3 + u_1 u_2 u_3 = 0.
\]
Under involutary linear mapping
\[
\left[\begin{array}{c}
x_0 \\ x_1 \\ x_2 \\ x_3
\end{array}\right] = 
\frac{1}{2}
\left[\begin{array}{cccc}
1 & 1 & 1 & 1 \\ 1 & 1 & -1 & -1 \\ 1 & -1 & 1 & -1 \\ 1 & -1 & -1 &  1
\end{array}\right]
\left[\begin{array}{c}
u_0 \\ u_1 \\ u_2 \\ u_3
\end{array}\right]
\]
the dehomogenized ($x_0=1$)
algebraic equation becomes
\[
p(x) = 
1-x_1^2-x_2^2-x_3^2-2x_1x_2x_3 =
\det \left[\begin{array}{ccc}1&x_1&x_2\\
x_1&1&x_3\\x_2&x_3&1\end{array}\right] =
\det F(x)
\]
which is the determinant of the 3x3 moment
matrix of the MAXCUT LMI relaxation.
The surface $p(x)=0$ is represented on Figure \ref{cayley_cubic},
using the {\tt surf} visualization package.
In particular, we can easily identify the convex connected
component containing the origin, described by
the LMI $F(x) \succeq 0$. The component has four vertices,
or singularities, for which the rank of $F(x)$ drops down to one.

In general, only curves and cubic surfaces admit
generically a determinantal representation. 
When $n>3$ or $m>3$ and no lifting is allowed,
the hypersurface $p(x)=0$ must be highly
singular to have a determinantal representation
\cite{beauville}, and hence, a fortiori, an LMI
representation.
This leaves however open the existence of
alternative algorithms consisting in constructing
symmetric linear determinantal representations
of modified polynomials $p(x)q(x)$, with
$q(x)$ globally nonnegative, say $q(x)=(\sum_i x^{2k}_i)$
or $(\sum_i x_i)^{2k}$
for $k\geq 1$ large enough.

Finally, let us conclude by remarking that, as a
by-product of the proof leading to Lemma \ref{rh},
checking numerically rigid convexity of
a scalar polynomial when $n>2$
amounts to checking positivity of a multivariate
Hermite matrix. See e.g. \cite{dumitrescu} for recent
developments on the use of semidefinite programming
for multivariate trigonometric polynomial matrix
positivity.

\section*{Acknowledgments}

This work benefited from feedback by (in alphabetical order)
Daniele Faenzi,
Leonid Gurvits, Fr\'ed\'eric Han,
Bill Helton, Toma\v z Ko\v sir, Salma Kuhlmann,
Jean-Bernard Lasserre, Adrian Lewis,
Bernard Mourrain, Pablo Parrilo,
Jens Piontkowski, Mihai Putinar, Bernd Sturmfels,
Jean Vall\`es and Victor Vinnikov.


\begin{thebibliography}{XX}

\bibitem{abhyankar}
S. S. Abhyankar, C. L. Bajaj.
Automatic parameterization of rational curves and surfaces III:
Algebraic plane curves. Computer Aided Geom. Design 5:309--321,
1988.

\bibitem{bpr}
S. Basu, R. Pollack, M.-F. Roy. Algorithms in 
real algebraic geometry. Springer, 2003.

\bibitem{beauville}
A. Beauville.
Determinantal hypersurfaces.
Michigan Math. J. 8:39--64, 2000.

\bibitem{becker}
B. Beckermann.
The condition number of real Vandermonde, Krylov and positive definite
Hankel matrices. Numer. Mathematik, 85:553-577, 2000.

\bibitem{bn}
A. Ben-Tal, A. Nemirovskii. Lectures on modern convex optimization.
SIAM, 2001.

\bibitem{multires}
L. Bus\'e, B. Mourrain, I.Z. Emiris, O. Ruatta, J. Canny, P. Pedersen,
I. Tonelli. Using the Maple {\tt multires} package.
INRIA Sophia-Antipolis, France, 2003.

\bibitem{kosir}
A. Buckley, T. Ko\v sir.
Determinantal representations of smooth cubic surfaces.
Geometriae Dedicata, 125:115--140, 2007.

\bibitem{callier}
F. M. Callier. On polynomial matrix spectral factorization
by symmetric extraction. IEEE Trans. Autom. Control,
30(5):453--464, 1985.

\bibitem{cox}
D. Cox, J. Little, D. O'Shea. Ideals, varietes, and
algorithms. Springer, 1992.

\bibitem{deconinck}
B. Deconinck, M. S. Patterson.
Computing the Abel map. Preprint, Dept. Applied Math.,
Univ. Washington, Seattle, 2007.

\bibitem{dvh}
B. Deconinck, M. van Hoeij.
Computing Riemann matrices of algebraic curves. Phys. D
152/153:123--152, 2001.

\bibitem{dixon}
A. C. Dixon. Note on the reduction of a ternary quantic to a
symmetrical determinant. Proc. Cambridge Phil.
Soc. 11:350--351, 1902.

\bibitem{dumitrescu}
B. Dumitrescu. Positive trigonometric polynomials and signal
processing applications. Springer, 2007.

\bibitem{edge}
W. L. Edge. Determinantal representations of $x^4+y^4+z^4$.
Proc. Cambridge Phil. Soc. 34:6--21, 1938.

\bibitem{mourrain}
M. Elkadi, B. Mourrain.
Introduction \`a la r\'esolution des syst\`emes polynomiaux.
Springer, 2007.

\bibitem{fk}
H. M. Farkas, I. Kra. Riemann surfaces. 2nd ed. Springer, 1992.

\bibitem{g02}
Y. Genin, Y. Hachez, Y. Nesterov, R. \c{S}tefan, P. Van Dooren, S. Xu.
Positivity and linear matrix inequalities. European J.
Control, 8(3):275--298, 2002.

\bibitem{g03}
Y. Genin, Y. Hachez, Y. Nesterov, P. Van Dooren.
Optimization problems over positive pseudo-polynomial matrices.
SIAM J. Matrix Anal. Appl., 25(1):57--79, 2003.

\bibitem{glr}
I. Gohberg, P. Lancaster, I. Rodman.
Matrix polynomials, Academic Press, 1982.

\bibitem{griffiths}
P. A. Griffiths. Introduction to algebraic curves.
AMS, 1989.

\bibitem{lfr}
S. Hecker, A. Varga, J. F. Magni.
Enhanced LFR toolbox for Matlab.
Proc. IEEE Symp. Computer Aided Control System Design,
Taiwan, 2004.

\bibitem{hmv}
J. W. Helton, S. A. McCullough, V. Vinnikov.
Noncommutative convexity arises from linear matrix inequalities.
J. Functional Analysis 240(1):105--191, 2006.

\bibitem{hv}
J. W. Helton, V. Vinnikov. Linear matrix inequality representation of
sets. Comm. Pure Applied Math. 60(5):654--674, 2007.

\bibitem{hn}
J. W. Helton, J. Nie.
Sufficient and necessary conditions for semidefinite representability
of convex hulls and sets.
{\tt arXiv:0709.4017}, September 2007.

\bibitem{hs}
D. Henrion, M. \v Sebek.
An efficient numerical method for the discrete-time symmetric matrix
polynomial equation. IEE Proc. Control Theory Appl. 145(5):443--448, 1998.

\bibitem{jk}
J. Je\v zek, V. Ku\v cera. Efficient algorithm for matrix spectral
factorization. Automatica, 21(6):663-669, 1985.

\bibitem{kaplan}
S. Kaplan, A. Shapiro, M. Teicher.
Several applications of B\'ezout matrices.
{\tt arXiv:math.AG/0601047}, January 2006.

\bibitem{k}
N. Kravistky. On the discriminant function of two commuting non-selfadjoint
operators. Integral Equ. Operator Theory, 3(1):97--124, 1980.

\bibitem{ks}
H. Kwakernaak, M. \v Sebek. Polynomial J-Spectral Factorization.
IEEE Trans. Autom. Control, 39(2):315--328, 1994.

\bibitem{lasserre}
J. B. Lasserre. Convex sets with lifted semidefinite representation.
LAAS-CNRS Research Report No. 07034, January 2007.

\bibitem{lpr}
A. S. Lewis, P. A. Parrilo, M. V. Ramana.
The Lax conjecture is true. Proc. Amer. Math. Soc. 133(9):2495--2499, 2005.

\bibitem{meyerbrandis}
T. Meyer-Brandis. Ber\"uhrungssysteme und symmetrische Darstellungen
ebener Kurven. Diplomarbeit, Univ. Mainz, 1998.

\bibitem{nn}
Y. Nesterov, A. Nemirovskii. Interior point polynomial algorithms in
convex programming, SIAM, 1994.

\bibitem{piontkowski}
J. Piontkowski.
Linear symmetric determinantal hypersurfaces.
Michigan Math. J. 54:117-146, 2006.

\bibitem{polyx}
PolyX, Ltd.
The Polynomial Toolbox for Matlab. 1998.

\bibitem{renegar}
J. Renegar.
Hyperbolic programs and their derivative relaxations.
Foundations Comput. Math. 6(1):59--79, 2006.

\bibitem{room}
T. G. Room.
The geometry of determinantal loci. Cambridge Univ. Press, 1938.

\bibitem{sendra}
J. R. Sendra, F. Winkler. Symbolic parametrization of curves.
J. Symbolic Comput. 12(6):607--631, 1991.

\bibitem{szilagyi}
I. Szil\'agyi.
Symbolic-numeric techniques for cubic surfaces. PhD Thesis,
RISC, Univ. Linz, Austria, 2005.

\bibitem{taussky}
O. Taussky. Nonsingular cubic curves as determinantal loci.
J. Math. Phys. Sci. 21(6):665--678, 1987.

\bibitem{tr}
H. L. Trentelman and P. Rapisarda.
New algorithms for polynomial J-spectral factorization.
Math. Control Signals Systems, 12:24--61, 1999.

\bibitem{tyrty}
E. E. Tyrtyshnikov. How bad are Hankel matrices ?
Numer. Math. 67:261--269, 1994.

\bibitem{vanhoeij}
M. van Hoeij. Rational parametrization of curves using canonical
divisors. J. Symbolic Comput. 23:209--227, 1997.

\bibitem{vinnikov86}
V. Vinnikov.
Self-adjoint determinantal representations of real irreducible
cubics. In: Operator Theory - Advances and Applications,
Vol. 19, Birkh\"auser, 1986.

\bibitem{vinnikov93}
V. Vinnikov.
Self-adjoint determinantal representations of real plane curves.
Math. Annalen, 296:453--479, 1993.

\bibitem{willems}
J. C. Willems.
Least squares stationary optimal control and the algebraic
Riccati equation. IEEE Trans. Autom. Control, 16(6):621--634, 1971.

\bibitem{yakubovich}
V. A. Yakubovich.
Factorization of symmetric matrix polynomials.
Dokl. Acad. Nauk. SSSR, 194(3):1261--1264, 1970.

\bibitem{zh}
J. C. Z\'u\~niga, D. Henrion.
A Toeplitz algorithm for polynomial J-spectral factorization.
Automatica, 42(7):1085--1093, 2006.

\end{thebibliography}
\end{document}